\documentclass[10pt]{article}
\usepackage{amsfonts}
\usepackage{amssymb}
\usepackage{amsmath}
\usepackage{amsthm}
\usepackage{graphicx}
\usepackage{tikz}

\usepackage{array,etoolbox}
\preto\tabular{\setcounter{magicrownumbers}{0}}
\newcounter{magicrownumbers}

\newtheorem{theorem}{Theorem}

\newtheorem{definition}[theorem]{Definition}

\newtheorem{lemma}[theorem]{Lemma}

\newtheorem{problem}{Problem}

\newtheorem{proof1}{Proof} 

\newenvironment{proofsketch}{{ \noindent \bfseries Proof Sketch}}{\qed} 

\newcommand{\itparen}[1]{\textup{(}#1\textup{)}}

\title{Notes on a short-cut to the proof of the $\mathbf{M}_3$-$\mathbf{N}_5$ Theorem}
\author{M. R. Emamy-K. \and Gustavo A. Mel\'{e}ndez R\'{i}os \\ Department of Mathematics, UPR R\'{i}o Piedras}
\date{February 22, 2024}

\begin{document}
\maketitle

\begin{abstract}
This paper presents two shortcuts to a classical proof of the $\mathbf{M}_3$-$\mathbf{N}_5$ Theorem, which can be found in \cite{DP} and \cite{BS}. To be precise, the shortcuts pertain a particular step of the proof that requires showing an algebraic equality. In addition, we  briefly discuss how to compare the lengths of the three proofs (the original and our two proposed shortcuts). To do so, we introduce two methods to compare the lengths of proofs based on algebraic lattice expressions. We call them the proof count method and the proof poset method. Both methods indicate that our proofs are shorter but the difference is more pronounced in the former. \\

\end{abstract}

\section{Introduction}

The $\mathbf{M}_3$-$\mathbf{N}_5$ Theorem is a basic theorem in lattice theory that characterizes all distributive and modular lattices. A classical proof of this theorem can be found in  \cite{DP} (also in \cite{BS}). Part of this proof consists of showing that a non-distributive modular lattice must have a sublattice isomorphic to $\mathbf{M}_3$. Showing this, in turn, requires proving an algebraic identity satisfied by three elements of the lattice. This paper proposes two shortcuts to the proof of this identity in the aforementioned references. They both consist of a reduction in the number of algebraic steps needed to show said identity. \\

We then use these shortcuts to discuss how to compare the lengths of proofs that consist of algebraic lattice expressions. After all, the length of a proof depends on how it is presented. To this end, we give an overview of results in \cite{msthesis} where we introduce two methods to compare the lengths of proofs:  the proof count method and the proof poset method. After applying them to the three proofs of the equation (the original in \cite{DP}, \cite{BS} and our two shortcuts), we get that our proofs are shorter but that the difference is more pronounced in the former method. \\

The rest of this paper is organized as follows. Section \ref{sec_back} provides some background material that includes a review of the $\mathbf{M}_3$-$\mathbf{N}_5$  Theorem with a sketch of its proof that shows where the algebraic equality of interest fits in it. We then present the original proof of the equality with our shortcuts in Section \ref{sec_3proofs}. Afterwards, we summarize our analysis of the lengths of the three proofs in Section \ref{sec_len_compare}. Finally, Section \ref{sec_concl} concludes with a summary of our results and ideas for continuing this work.

\section{Background} \label{sec_back}

We set the stage for our proposed proof shortcuts in two steps. First, we briefly review some lattice theory basics and establish the notation we use in Subsection \ref{ssec_prelim}. Afterwards, Subsection \ref{ssec_m3n5} reviews the $\mathbf{M}_3$-$\mathbf{N}_5$  Theorem and gives a sketch of the proof in \cite{DP} and \cite{BS}. For more information regarding lattices and posets, we refer the reader to \cite{DP}.

\subsection{Preliminaries} \label{ssec_prelim}

We begin with some notation. Given a lattice $L$, we denote its meet operation with $\land$ and its join operation with $\lor$. In addition, we denote its induced order with $\leq$, using $<$ for strict order. If two lattices $L$ and $K$ are isomorphic, we write $L \cong K$. Finally, recall that $\mathbf{M}_3$ and $\mathbf{N}_5$ denote the diamond and pentagon lattices, respectively. These are shown in Figures \ref{m3} and \ref{n5}. \\

\begin{figure}
\begin{minipage}{0.45\textwidth}
\centering
\begin{tikzpicture}
[roundnode/.style={circle, draw=black, fill=white, very thin, text width =4.00mm, inner sep = 0pt, text centered}]
 \node [roundnode] (1) at (0,2) {$1$};
  \node [roundnode] (a) at (-1,1) {$a$};
  \node [roundnode] (b) at (0,1) {$b$};
  \node [roundnode] (c) at (1,1) {$c$};
  \node [roundnode] (0) at (0,0) {$0$};
 \draw (1) -- (a) -- (0) -- (b) -- (1) -- (c) -- (0);
\end{tikzpicture}
\caption{The diamond $\mathbf{M}_3$} \label{m3}
\end{minipage} \hfill
\begin{minipage}{0.45\textwidth}
\centering
\begin{tikzpicture}
[roundnode/.style={circle, draw=black, fill=white, very thin, text width =4.00mm, inner sep = 0pt, text centered}]
 \node [roundnode] (1) at (0,2) {$1$};
  \node [roundnode] (u) at (-1,1.4) {$u$};
  \node [roundnode] (v) at (-1,.6) {$v$};
  \node [roundnode] (w) at (1,1) {$w$};
  \node [roundnode] (0) at (0,0) {$0$};
 \draw (1) -- (u) -- (v) -- (0) -- (w) -- (1);
\end{tikzpicture}
\caption{The pentagon $\mathbf{N}_5$} \label{n5}
\end{minipage}
\end{figure}

We also recall some basic properties and results about lattices. First, the lattice operations $\land$ and $\lor$ satisfy four axioms: associativity, commutativity, idempotency, and absorption. Second, all distributive lattices are modular. Last, we formally state two results that we will make reference to.

\begin{lemma} \itparen{Connecting Lemma} \label{connect-lemma}
Let $L$ be a lattice with $a,b \in L$ and induced order $\leq$. Then the following are equivalent:
\begin{enumerate}
\item $a \leq b$;
\item $a \lor b = b$;
\item $a \land b = a$.
\end{enumerate}
\end{lemma}

\begin{lemma} \itparen{half-way lemma} \label{half-lemma}
Let $L$ be a lattice with $a,b,c \in L$. Then
\begin{enumerate}
\item $a \land (b \lor c) \geq (a \land b) \lor (a \land c);$
\item $a \geq c \implies a \land (b \lor c) \geq (a \land b) \lor c.$
\end{enumerate}
\end{lemma}

\subsection{$\mathbf{M}_3$-$\mathbf{N}_5$ Theorem} \label{ssec_m3n5}

We now review the $\mathbf{M}_3$-$\mathbf{N}_5$ Theorem and sketch the proof in \cite{DP} and \cite{BS}. This theorem characterizes distributive and modular lattices. The idea is that $\mathbf{N}_5$ is a forbidden sublattice of modular lattices in the sense that all lattices that contain it are non-modular. Similarly, $\mathbf{M}_3$ and $\mathbf{N}_5$ are forbidden sublattices of distributive lattices.

\begin{theorem} \itparen{$\mathbf{M}_3$-$\mathbf{N}_5$ Theorem} \label{m3n5} 
Let $L$ be a lattice.
\begin{enumerate}
\item $L$ is non-modular if and only if $\mathbf{N}_5$ is a sublattice of $L$.
\item $L$ is non-distributive if and only if $\mathbf{N}_5$ or $\mathbf{M}_3$ is a sublattice of $L$.
\end{enumerate}

\begin{proofsketch}
We provide a sketch of the proof in Davey and Priestley \itparen{page 89 of \cite{DP}; also in \cite{BS}}. Given that sublattices of modular \itparen{distributive} lattices are modular \itparen{distributive} and that $\mathbf{N}_5$ and $\mathbf{M}_3$ are non-modular and modular but not distributive respectively, it suffices to show that \itparen{i} a non-modular lattice has $\mathbf{N}_5$ as a sublattice and \itparen{ii} that a non-distributive modular lattice has $\mathbf{M}_3$ as a sublattice.\\

\indent Part \itparen{i}: Suppose that L is a non-modular lattice. Then there exist $d, e, f$ in $L$ with $d > f$ such that $d \land (e \lor f) > (d \land e) \lor f$. This is because by Lemma \ref{half-lemma},
\begin{equation}
d \land (e \lor f) \neq (d \land e) \lor f \implies d \land (e \lor f) > (d \land e) \lor f.
\end{equation}
Define $a = (d \land e) \lor f$ and $b =  d \land (e \lor f)$. It can be shown that $a \land e = b \land e$ and $a \lor e = b \lor e$. Let $p = a \land e$ and $q = a \lor e$. Then $K = \{p,q,a,b,e\}$ is a sublattice of $L$ that satisfies the necessary joins and meets to be isomorphic to $\mathbf{N}_5$. It can be verified that all 5 elements of $K$ are distinct. Therefore, $K \cong \mathbf{N}_5$ and $L$ has a sublattice isomorphic to $\mathbf{N}_5$.\\

\indent Part \itparen{ii}: Suppose that $L$ is a modular but non-distributive lattice. Then there exist $d, e, f \in L$ such that $d \land (e \lor f) > (d \land e) \lor (d \land f)$. This is because by Lemma \ref{half-lemma},
\begin{equation}
d \land (e \lor f) \neq (d \land e) \lor (d \land f) \implies d \land (e \lor f) > (d \land e) \lor (d \land f).
\end{equation}
Define the following five elements in $L$:
\begin{align}
p &= (d \land e) \lor (e \land f) \lor (f \land d),  \\
q &= (d \lor e) \land (e \lor f) \land (f \lor d), \\
u &= (d \land q) \lor p, \\
v &= (e \land q) \lor p, \\
w &= (f \land q) \lor p.
\end{align}
It can be shown that $K = \{p,q,u,v,w\}$ is a sublattice of $L$ isomorphic to $\mathbf{M}_3$. This requires proving that
\begin{enumerate}
\item $p < q;$
\item $p \leq u,v,w \leq q;$
\item $u \land v = v \land w = w \land u = p;$
\item $u \lor v = v \lor w = w \lor u = q;$
\item All five elements of K are distinct.
\end{enumerate}
This concludes our proof sketch.
\end{proofsketch}
\end{theorem}

We have finally arrived at our point of interest. Note that item 3 in the list above requires showing three identities:  $u \land v = p$,  $v \land w = p$, and $w \land u = p$. We remark that the proofs of all three of them are basically the same. Therefore, we focus on the proof of the identity $u \land v = p$ for which we propose the two shortcuts mentioned at the beginning of this paper. Thus, we have established the context of our shortcuts within the proof of the $\mathbf{M}_3$-$\mathbf{N}_5$ Theorem and are ready to present them.

\section {3 Proofs of $u \land v = p$} \label{sec_3proofs}

We now discuss the three proofs of $u \land v = p$, the original in \cite{DP} and \cite{BS} plus our two shortcuts. All three have the same general strategy: algebraic manipulation of lattice expressions starting from $u \land v$ and ending with $p$. In particular, they involve substitution of $u$, $v$, $p$, and $q$ for their definitions and applications of lattice axioms, modular law, and the Connecting Lemma. \\

\indent Each proof is followed by a list of justifications for each algebraic step. In these justifications we make several references to the modular law, that is; for all $a,b,c \in L$ with $a \geq c,$ we have   $a \land (b \lor c) = (a \land b) \lor c$ or without distinction to its dual, that is; for all $a \leq c,$  we have   $a \lor (b \land c) = (a \lor b) \land c$ (we also underline $a$ and $c$ just before each step when modularity is applied to facilitate  the reading).  Not all uses of commutativity and associativity in the proofs are mentioned. The reader is advised to consult the proof sketch of Theorem \ref{m3n5}  to recall the definitions of and the relations among the elements treated in the proofs. \\

\indent We begin with the classical proof given in page 91 of \cite{DP} (book by B. A. Davey and H. A. Priestley). It is also in \cite{BS}. We will call it Proof 1.

\begin{proof1} \itparen{Classical proof in \cite{DP} and \cite{BS}}
\begin{subequations}
\begin{align}
u \land v &= ((d \land q) \lor \underline{p}) \land (\underline{(e \land q) \lor p}) \label{dp1} \\
&= (((e \land \underline{q}) \lor \underline{p}) \land (d \land q)) \lor p  \label{dp2} \\
&= ((q \land (e \lor p)) \land (d \land q)) \lor p  \label{dp3} \\
&= ((e \lor p) \land (d \land q)) \lor p  \label{dp4} \\
&= ((d \land (e \lor f)) \land (e \lor (f \land d))) \lor p  \label{dp5}\\
&= (d \land ((\underline{e \lor f}) \land (\underline{e} \lor (f \land d)))) \lor p  \label{dp6}\\
&= (d \land (((e \lor f) \land (f \land d)) \lor e )) \lor p  \label{dp7}\\
&= (\underline{d} \land ((\underline{f \land d}) \lor e )) \lor p  \label{dp8}\\
&= ((d \land e) \lor (f \land d)) \lor p = p \label{dp9}
\end{align}
\end{subequations}

\underline{Justification for each step of Proof 1}

\begin{itemize}
\item \itparen{\ref{dp1}} By definition of $u$ and $v$.
\item \itparen{\ref{dp2}} Applying modular law with $a = (e \land q) \lor p$, $b = d \land q$, and $c = p$.
\item \itparen{\ref{dp3}} Applying modular law with $a = q$, $b = e$, and $c = p$.
\item \itparen{\ref{dp4}} By associativity, commutativity, and idempotency on $q$.
\item \itparen{\ref{dp5}} By definitions of $p$ and $q$, associativity, commutativity, and absorption.
\item \itparen{\ref{dp6}} By associativity.
\item \itparen{\ref{dp7}} Applying modular law with $a = e \lor f$, $b = f \land d$, and $c = e$.
\item \itparen{\ref{dp8}} Since $f \land d \leq f \leq e \lor f$.
\item \itparen{\ref{dp9}} Applying modular law with $a = d$, $b = e$, and $c = f \land d$ and using the fact that $(d \land e) \lor (f \land d) \leq p$ by definition of $p$. \qed
\end{itemize}
\end{proof1}

We now present our two new proofs of the same result. These will be Proofs 2 and 3 respectively.  We remark that the versions of both presented here are an improvement from those in \cite{msthesis}.
\begin{enumerate}
\item Proof 2: It takes a short-cut in Proof 1.
\item Proof 3: It results from trying to re-create either of the first two proofs from memory but taking a different route to the same result instead.
\end{enumerate}

\begin{proof1} 
\begin{subequations}
\begin{align}
u \land v &= ((d \land q) \lor \underline{p}) \land (\underline{(e \land q) \lor p}) \label{prof1}\\
&= (((e \land \underline{q}) \lor \underline{p}) \land (d \land q)) \lor p \label{prof2}\\
&= ((q \land (e \lor p)) \land (d \land q)) \lor p \label{prof3}\\
&= ((e \lor p) \land (d \land q)) \lor p \label{prof4}\\
&= ((\underline{d} \land (e \lor f)) \land (e \lor (\underline{f \land d}))) \lor p \label{prof5}\\
&= ((e \lor f) \land ((d \land e) \lor (f \land d))) \lor p \label{prof7} = p
\end{align}
\end{subequations}

\underline{Justification for each step of Proof 2}

\begin{itemize}
\item \itparen{\ref{prof1}}-\itparen{\ref{prof5}} Identical to \itparen{\ref{dp1}}-\itparen{\ref{dp5}} in Proof 1.
\item \itparen{\ref{prof7}} By commutativity, associativity, applying modular law with $a = d$, $b = e$, and $c = f \land d$, and the fact that $(e \lor f) \land ((d \land e) \lor (f \land d)) \leq p$ by definition of $p$. \qed
\end{itemize}
\end{proof1}

\begin{proof1} 
\begin{subequations}
\begin{align}
u \land v &= [(d \land \underline{q}) \lor \underline{p}] \land [(e \land \underline{q}) \lor \underline{p}] \label{def-uv}\\
&= [q \land (d \lor p)] \land [q \land (e \lor p)] \label{modx2}\\
&= q \land \{ (d \lor p) \land (e \lor p) \} \label{com-assoc-idem}\\
&= q \land \{ [\underline{d \lor (e \land f)}] \land [e \lor (\underline{d \land f})]\} \label{p-def+abs}\\
&= q \land \{[(d \lor (\underline{e \land f})) \land \underline{e}] \lor (d \land f)\} \label{mod2}\\
&= q \land \{ [(e \land d) \lor (e \land f)] \lor (d \land f)\} \label{mod3}\\
&= q \land \{ (e \land d) \lor (e \land f) \lor (d \land f)\} = p \label{assoc}
\end{align}
\end{subequations}

\underline{Justification for each step of Proof 3}

\begin{itemize}
\item \itparen{\ref{def-uv}} By definition of $u$ and $v$.
\item \itparen{\ref{modx2}} Applying modular law twice:
\begin{enumerate}
\item With $a = q$, $b = d$, and $c = p$.
\item With $a = q$, $b = e$, and $c = p$.
\end{enumerate}
\item \itparen{\ref{com-assoc-idem}} By commutativity, associativity, and idempotency.
\item \itparen{\ref{p-def+abs}} By definition of $p$ and absorption.
\item \itparen{\ref{mod2}} Applying modular law with $a = d \lor (e \land f)$, $b=e$, and $c = d \land f$. Modular law can be applied because $a = d \lor (e \land f)\geq d \geq d \land f = c$.
\item \itparen{\ref{mod3}} Applying modular law with $a = e $, $b=d$, and $c = e \land f$.
\item \itparen{\ref{assoc}} By associativity, the definition of $p$, and $p < q$. \qed
\end{itemize}
\end{proof1}

\section{Comparing Proof Lengths} \label{sec_len_compare}

Having shown the three proofs of the identity $u \land v = p$, we turn our attention to comparing their length. It can readily be seen that Proofs 1, 2, 3 respectively have 9, 6, and 7 lines. However, we know that the length of a proof is not an immutable property since it depends on the level of detail given (and hence on presentation). Thus, we study this matter further by coming up with two methods to compare the length of the three proofs. We will call them proof count and proof poset respectively. The first will be na\"{i}ve while the second will be more rigorous. Both methods are of our own design. Here, we limit ourselves to an overview of the results of our analysis. The details can be found in \cite{msthesis}.
\subsection{Proof Count Method}
\indent The proof count method is a straightforward brute-force approach. We measure the length of each proof by counting the lines and characters it uses. To be more precise, we count the following:
\begin{enumerate}
\item lines,
\item variables,
\item operation symbols,
\item grouping symbols,
\item equal signs,
\item all symbols.
\end{enumerate}

Table \ref{proof_count} shows the results of comparing the three proofs with the proof count method. The main observation we make is that Proofs 2 and 3 use significantly less symbols than Proof 1 (129 and 155 vs. 185). There are two discrepancies between the data here and that in \cite{msthesis} that must be explained. The first is that we moved the final $=p$ of each proof to the previous line for the sake of efficiency. The second is that we made improvements to Proofs 2 and 3, as mentioned in Section \ref{sec_3proofs}.

\begin{table}
\centering
\begin{tabular}{| c | c | c | c | c | c | c |}
\hline
Proof & lines & variables & operations & grouping & equal & \textbf{symbol total}\\
\hline
Proof 1 & 9 & 57 & 46 & 72 & 10 & \textbf{185} \\
\hline
Proof 2 & 6 & 40 & 32 & 50 & 7 & \textbf{129} \\
\hline
Proof 3 & 7 & 48 & 39 & 60 & 8 & \textbf{155} \\
\hline
\end{tabular}
\caption{Results of Proof Count}  \label{proof_count}
\end{table}

\subsection{Proof Poset Method}

The second method we use is the proof poset method. In it, we build a poset (specifically a chain) to represent each of the proofs as follows:

\begin{enumerate}
\item Top element: $u \land v = [(d \land q) \lor p] \land [(e \land q) \lor p]$;
\item Bottom element: $p$;
\item Each vertex is a statement in proof;
\item Vertex $i$ covers vertex $j$ if we can go from statement $i$ to statement $j$ by applying only one rule from a small list of basic rules (to be given in Definition \ref{rule_list});
\item Must add statements to proof if going from $i$ to $i+1$ requires more than one rule.
\end{enumerate}

The idea is to examine in detail the logical structure of each of the proofs by decomposing them into their building blocks. We admit the following list of basic rules for the covering relation of the proof posets.

\begin{definition} \itparen{basic rules for proof poset} \label{rule_list}
\begin{itemize}
\item Lattice axioms
\begin{enumerate}
\item L1: associative laws
\item L2: commutative laws
\item L3: idempotency laws
\item L4: absorption laws
\end{enumerate}
\item Modular law
\item Other rules
\begin{enumerate}
\item Def: Substitute an element for its definition or vice-versa.
\item Applying Connecting Lemma or some other order property of lattices.
\end{enumerate}
\end{itemize}
\end{definition}

In order to shorten the posets, we will also allow the following combinations of the modular law with commutativity as a single step (covering relation of poset).

\begin{definition} \itparen{list of combinations of modular law with commutativity} \label{mod_comm_comb}
\begin{enumerate}
\item M: Modular law
\begin{equation}
a \geq c \implies a \land (b \lor c) = (a \land b) \lor c
\end{equation}
\item M1: Modular law with inner commute
\begin{equation}
a \geq c \implies a \land (c \lor b) = (a \land b) \lor c
\end{equation}
\item M2: Modular law with outer commute
\begin{equation}
a \geq c \implies (b \lor c) \land a = (a \land b) \lor c
\end{equation}
\item M3: Modular law with reverse inner commute
\begin{equation}
a \geq c \implies (b \land a) \lor c = a \land (b \lor c)
\end{equation}
\item M4: Modular law with reverse outer commute
\begin{equation}
a \geq c \implies c \lor (a \land b) = a \land (b \lor c)
\end{equation}
\end{enumerate}
\end{definition}

Needless to say, none of the proofs given above is fully decomposed based on the lists of Definitions \ref{rule_list} and \ref{mod_comm_comb}. Hence, they require adding intermediate statements. Doing so results in having three fully decomposed proofs from which we can construct three proof posets as indicated above. The fully decomposed proofs and the proof posets can be found in  \cite{msthesis}. \\

We now give an overview of the results from the proof poset method. We compute the lengths of the decomposed proofs (number of vertices of poset) in Table \ref{proof_poset} and count how many times each proof uses each basic rule in Table \ref{rule_count}. From Table \ref{proof_poset}, we get that Proofs 2 and 3 are only slightly shorter than Proof 1 based on our notion of ``basic rule"  in Definitions \ref{rule_list} and \ref{mod_comm_comb}. The main observation from Table \ref{rule_count} is that Proof 2 uses the modular law one time less than the other proofs (3 vs. 4). Noting that the modular law is less basic than the axioms, this could be interpreted as saying that Proof 2 is somehow logically simpler than the other two. \\

\indent Therefore, we conclude the following from our proof comparison exercise:
\begin{enumerate}
\item Proof Count Method: Proofs 2 and 3 are shorter than Proof 1.
\item Proof Poset Method: Proofs 2 and 3 are slightly shorter than Proof 1.
\item Proof 2 is the shortest proof either way.
\end{enumerate}

\begin{table}
\centering
\begin{tabular}{| c | c |}
\hline
Proof 1 & 32 \\
\hline
Proof 2 & 29 \\
\hline
Proof 3 & 30 \\
\hline
\end{tabular}
\caption{Length of proof posets} \label{proof_poset}
\end{table}

\begin{table}
\centering
\begin{tabular}{| c | c | c | c | c | c | c | c | c | c | c | c | }
\hline
   Proof & L1 & L2 & L3 & L4 & M & M1 & M2 & M3 & M4 & Def & Other \\
     \hline
   1 &  9 & 9 & 1 & 4 & 0 & 2 & 1 & 1 & 0 & 2 & 2 \\
     \hline
   2 &  9 & 8 & 1 & 4 & 1 & 0 & 1 & 1 & 0 & 2 & 1 \\
     \hline
   3 &  10 & 6 & 1 & 4 & 1 & 0 & 1 & 2 & 0 & 3 & 1 \\
    \hline
\end{tabular}
\caption{Count of basic rules invoked by Proofs 1-3} \label{rule_count}
\end{table}

\section{Conclusion} \label{sec_concl}

We bring this paper to a close with a summary of our results and some suggestions for continuing this work. We presented two shortcuts to a part of a classical proof of the $\mathbf{M}_3$-$\mathbf{N}_5$ Theorem found in \cite{DP,BS}. These were the result of ``algebraic efficiencies" in proving the equation $u \land v = p$ in said proof. We also gave an overview of our work comparing the length of the three proofs in \cite{msthesis} where we proposed and applied two methods for doing so: the proof count method and the proof poset method. Both demonstrated that the new proofs were shorter than the proof in \cite{DP,BS} but the difference was more marked in the proof count method. \\

Needless to say, it is possible to keep pushing the boundaries of this work. In particular, we can develop further both proof length comparison methods and apply them to other proofs in lattice theory. This could lead to a better-defined notion of ``proof length" and, thus, a potential classification of lattice theory proofs based on it. We finish with some ideas for the development of both methods.

\begin{problem} \itparen{Proof count method ideas}
\begin{enumerate}
\item Consider other aspects that can be counted.
\item Add different weights to different symbols.
\end{enumerate}
\end{problem}

\begin{problem} \itparen{Proof poset method ideas}
\begin{enumerate}
\item Refine better basic rules. For instance, break the combined modular law rules (M1-M4) to see the impact on the final results.
\item Add different weights to different rules based on how basic they are.
\end{enumerate}
\end{problem}




\end{document}